\documentclass[a4paper, 12pt] {article}

\usepackage{fende}  

\begin{document}

\title{RCF\,4 \\
       Inconsistent Quantification \\ 
       $\mr{AC}\ \forall\exists!\ \eps$}

         \footnotetext{
           this is part 4 
           of a cycle on \emph{Recursive Categorical Foundations},
           we rely on RCF\,1 and quote section 2 of RCF\,2}
           
         \footnotetext{legend of \textsc{logo:} 
           $\AC:$ Axiom of Choice,\ $\forall\exists!:$ 
           Discrete map definition by 
           ``\,$\forall_a\,\exists!_b\,\ph(a,b)$\,'', \\
           $\eps:$ Iterative evaluation of $\PRa$ map codes, 
           $\PRa$ the Theory of Primitive Recursion with
           predicate abstraction}

    \author{Michael Pfender} 

    \date{version 2\footnote{
                      update to version 1 mainly 
                      of section 5 on \emph{self-evaluation}},
                     \ September 2009\,\footnote{last revised \today}}     

\maketitle

\abstract
{We exhibit canonical \emph{middle-inverse} Choice maps within 
categorical (Free-Variable) Theory of Primitive Recursion
as well as in Theory of \emph{partial} PR maps over
Theory  of Primitive Recursion with \emph{predicate abstraction.} 
Using these choice-maps, defined by $\mu$-recursion, we address the 
\emph{consistency problem} for a minimal Quantified extension 
$\bfQ$ of latter two theories: We prove, that $\bfQ$'s 
$\exists$-defined $\mu$-operator coincides on PR predicates with 
that inherited from theory of \emph{partial} PR maps.   
We strengthen Theory $\bfQ$ by axiomatically forcing the
lexicographical order on its $\omega^{\omega}$ to become a well-order:
``finite descent''. Resulting theory admits 
\emph{non-infinit PR-iterative descent} schema $(\pi)$ which constitutes 
Cartesian PR Theory $\piR$ introduced in RCF\,2. 

A suitable Cartesian subSystem of $\bfQ+\wo(\omega^{\omega})$
above, extension of $\piR$ ``inside'' Theory $\bfQ+\wo(\omega^{\omega}),$ 
is shown to admit code \emph{self-evaluation:} extension of formally 
\emph{partial} code evaluation of $\piR.$ Appropriate 
\emph{diagonal argument} then shows inconsistency of this 
subSystem and (hence) of its extensions $\bfQ+\wo(\omega^{\omega})$ 
and $\ZF.$}

\section{Introduction}  

We begin with Proof of a local, \emph{middle-inverse} form $\ACCmi$
of---Countable---Choice. This for \emph{fundamental} Free-Variables 
(categorical) Theory $\PR,$ as well as for Theory 
$\hatPRa = \widehat{\PRa}$ of \emph{partial} maps over Theory \\ 
$\PRa = \PR+(\abstr)$ of Primitive Recursion with 
\emph{predicate abstraction} 
$\bfan{\chi: A \to 2} \bs{\mapsto}\,\emph{Object}\ \set{A\,|\,\chi}.$
Equational (!) Axiom $\ACCmi$ is preserved by theory strengthening,
and by theory \textbf{extension}---the latter with respect to 
$\PRa$-defined maps.

\smallskip
$[\,\AC$ cannot hold for Theory $\PRa$ itself consistently$\,]$
  
\smallskip
What we \emph{can} prove is ``even'' middle inverse form 
$\ACCmi$ of $\AC,$ for ``classically'' quantified Arithmetical 
Theory $\bfQ = \PRa+\forall\exists!,$ having (possibility of) 
\emph{``discrete'' map-definition,} via left-total, right-unique 
binary predicates $\ph = \ph(a,b): A \times B \to 2,$ a possibility
for map-definition typical for \textbf{set} theorie(s). 
 
\smallskip
For Ordinal $\N[\omega] \subset \N^* = \omega^{\omega}$ we recall 
schema $(\pi) = (\pi_{\N[\omega]})$ of \emph{finite descent} for 
\emph{Complexity Controlled Iteration} with \emph{complexity values}
in $\N[\omega],$ and definition of strengthening $\piR = \PRa+(\pi)$ 
of $\PRa:$ within $\piR$ the defined-arguments enumerations of
its $\mrCCI$'s are forced to become epi, ``onto'': These
$\mrCCI$'s \emph{on-terminate} within Theory $\piR,$
in particular so does formally 
\emph{partial}---iterative---\emph{code evaluation} of 
Theory $\piR,$ \cf part RCF\,2.

\smallskip
``Critical'' Theory, namely Theory 
$\bfQ^{\wo}:$ $\hatPRa \bs{\sqsup} \PRa,$
enriched by existential Quantification giving ``total'' predicates
$\exists_n\ph(a,n): A \to 2$ from ``total''---$\PRa$---predicates 
$\ph: A \times \N \to 2,$ makes the (canonical) 
middle-inverse \emph{partial} maps, middle-inverse to \emph{defined-arguments} 
enumerations of $\piR$'s $\mrCCI$'s, into ``total'' maps, \emph{maps} 
within $\bfQ^{\wo}.$ Adding these $\bfQ^{\wo}$-maps as ``\emph{total}''
maps to Theory $\piR,$ \ie forcing by $\bfQ^{\wo}$-consistent 
\textbf{axiom} the enriched Theory---a priori only PR 
\emph{monoidal}---to become \emph{Cartesian,} allows for resulting 
(Cartesian PR) Theory $\dotpiR$ code \emph{self-evaluation}
$\doteps(u,a): \mrdotpiR \times \X \to \X$ (within $\dotpiR$).

\medskip
From this then results---by appropriate 
\emph{diagonal argument}---inconsistency of $\dotpiR$ as well as of
its extensions $\bfQ^{\wo},\ \PA+\wo(\omega^{\omega})$ and $\ZF.$

\section{Middle-Inverse Choice Maps in Theories $\PR$ and $\hatPRa$}

\textbf{Definition:} 
For a given map (term) $f = f(a): A \to B$ of a (categorical) theory 
$\T,$ a $\T$-map $f': B \to A$ is called a \emph{Choice map} 
for $f,$ in the \emph{middle-inverse} sense, if  
  $$f = f \,\circ\, f' \,\circ\, f: A \xto{f} B \xto{f'} A \xto{f} B
                                                      \ \text{in}\ \T.$$
If the given $\T$ map $f: A \onto B$ is a $\T$-\emph{epi,} then obviously 
$f': B \into A$ is a $\T$-\emph{section} for $f.$  

\smallskip
\textbf{Definition:} A (categorical) Theory $\T$ with \emph{terminal} 
Object $\one$---or at least a \emph{half-terminal Object} $\one:$ 
each Object $A$ admits a (non-necessary unique) $\T$-map 
$!: A \to \one$---is said to admit \emph{(middle-inverse) Choice,} or to satisfy 
Axiom $\ACmi$ if each $\T$ map $f: A \to B$ coming with a 
\emph{point} $a_0: \one \to A,$ admits a \emph{middle inverse} 
map $f': B \to A$ in the sense above.

 \newpage

\medskip
\textbf{Remarks:}
\begin{itemize}
\item
If $\T$ satisfies $\ACmi,$ then each ``pointed'' $\T$-\emph{epi} 
is a \emph{retraction:} $\T$ satisfies the (local) 
\emph{Axiom of Choice} $\AC.$ 
And---dually---each pointed $\T$-\emph{mono} then is a \emph{section.}

\item
In \textbf{set} theories, requirement of \emph{pointed Domains} seems
to be redundant, since \emph{non-empty} sets \emph{have} points, by
extensionality axiom. But are these points \emph{available} for 
``construction'' below, without (set-theoretical) Axiom $\ACC$
of \emph{Countable Choice}?
In our case \emph{yes:} by the ``set-theoretical'' $\mu$-operator,
available \eg in $\PA:$ ``(Classical) $\PA = \PRA+\exists$\,''



\end{itemize}

\bigskip
\textbf{Countable Choice Theorem for $\PR$ and $\hatPRa:$}
\begin{enumerate} [(i)]

\item 
Fundamental theory $\PR$ of Primitive Recursion---Objects: finite 
(binary bracketed) powers of $\N,$ not yet formal extensions 
(abstractions) $\set{A\,|\,\chi: A \to 2}$ --, admits, 
within itself, middle-inverse \emph{Choice maps} $f': B \to A$ 
for all of its maps $f: A \to B.$  

In particular, all epis of this fundamental theory turn out to be 
retractions: $\PR$ satisfies $\AC$ (here $\ACC$).

\smallskip
$[\,$All Objects $A$ of $\PR$ are \emph{pointed,} 
by---componentwise defined---zero $0: \one \to A.$ We just need
any point. $\PR$ is not a \emph{pointed category,} since maps are not
required to map ``canonical'' points into canonical ones.$]$

\item 
Theory $\hatPRa,$ of \emph{partial} PR maps over \emph{basic} Theory
$\PRa = \PR+(\abstr)$ of Primitive Recursion with 
\emph{predicate abstraction,} again admits axiom $\ACC$ of (Countable)
Choice, in the form of middle-inverse \emph{partial} PR maps to
arbitrary partial PR maps.


\item
Middle-inverse form $\ACmi$ of $\AC$ is clearly inherited by 
\emph{strengthenings} of a theory, because of its purely equational 
character: To each map is associated a map in the converse direction, 
with ``characteristic'' middle-inverse \emph{equation}---maintained. 

\item
\textbf{Problem:} Does Theory $\PRa$ ``itself'' admit $\AC?$

Middle-inverse $f^-: B \parto A$ to a $\PRa$ map $f: A \to B$ is 
in general not $\PRa.$ Use of $\AC$ in its \emph{epis-have-sections}
form cannot be inherited by $\PRa$ from $\hatPRa$ since $\PRa$
epis are a priori not $\hatPRa$ epis: a ``direct'' proof would need
$\PRa$ $\pb$'s to pull back epis into epis, and this is
excluded in general, by an argument discussed in part RCF\,2.

\end{enumerate}
 
\textbf{Proof} of assertion (i) by recursive case distinction on the 
structure of $f: A \to B$ in $\PR,$ $A,B$ \emph{fundamental,} \ie
of form of a (binary bracketed) finite power of object $\N:$
\begin{itemize}

\item Case of map-constants: All of these come with retractions
or with sections, in particular since each of the fundamental (!)
Objects $A$ comes with a (componentwise defined) zero 
$0_A: \one \to A.$










\item Composition $f = h \circ g: A \to B \to C:$
$f' \defeq g' \circ h': C \to B \to A.$

\item Cylindrification $f = (C \times g): C \times A \to C \times B:$ \\
$f' \defeq (C \times g'): C \times B \to C \times A.$ 

\item Iteration $f = g^\S = g^\S(a,n): A \times \N \to A:$ \\
$(\id_A,0 \circ !_A): A \to A \times \N$ is a section to $f.$

\end{itemize}

\textbf{Proof} of assertion (ii): \emph{middle-inverse Choice} 
for Theory $\hatPRa:$

For $f = \an{(d_f,\widehat{f}): D_f \to A \times B}: A \parto B$
within $\PRa,$ we could choose middle-inverse just (graph-) opposite 
to $f,$ namely 
  $$f^- \bydefeq 
      \an{(\widehat{f},d_f): D_f \to B \times A}: B \parto A.$$
But wanted \textbf{proof} of middle-innverse property
  $$f \parcirc g \parcirc f = f \parcirc f^- \parcirc f \pareq f:
                                    A \parto B \parto A \parto B$$
is more conceptual---and simpler---if we use definition of
partial maps inside $\PRa$ via $\mu$-recursion, \cf RCF\,1:

We \textbf{define} our middle-inverse candidate 
$g = g(b): B \parto A$ in $\muR \bs{\iso} \hatPRa$ as follows, 
(essentially) via a (partial) $\muR$-map 
  $$\mu_g = \mu_g(b) \defeq 
      \mu\set{\hat{a} \in D_f\,|\,\widehat{f} (\hat{a}) \doteq_B b}:
           B \parto D_f,$$
this with respect to canonical, \NAME{Cantor} ordering of Object
$D_f = \set{D\,|\,\zeta}$ inherited from $\N$ via $D$ fundamental.

Partial map $g: B \parto A$ is then choosen as 
\begin{align*}
& g \defeq d_f \parcirc \mu_g: B \parto D_f \to A,
                                  \ \text{with}\ b \in B\ \free: \\
& g(b) \defeq d_f(\mu\set{\hat{a} \in D_f\,
                            |\,\widehat{f} (\hat{a}) \doteq_B b}):
                                                          B \parto A.
\end{align*}
This $g: B \parto A$ \emph{is} a middle-inverse to $f: A \parto B,$ 
since---preliminary result:
\begin{align*}
& f \parcirc g \parcirc f \parcirc d_f 
        = f \parcirc g \parcirc f \parcirc d_f(\hat{a}) \\
& \pareq f \parcirc g \parcirc \widehat{f} (\hat{a}) \\
& = f \parcirc d_f \parcirc
      (\mu\set{\hat{a}' \in D_f\,|\,\widehat{f} (\hat{a}') 
                                      \doteq_B \widehat{f} (\hat{a})}) \\
& = \widehat{f} \circ 
      (\min\set{\hat{a}' \leq_{D_f}\,\hat{a}\,|
                 \,\widehat{f} (\hat{a}') \doteq_B \widehat{f} (\hat{a})}) \\
& = \widehat{f} (\hat{a}) \pareq f \parcirc d_f(\hat{a})
           = f \parcirc d_f: D_f \to A \parto B.
\end{align*}
In order to get rid of the leading $d_f: D_f \to A$ on both sides
of the (resulting) $\hatPRa$ equation above, we use the commuting
$\hatPRa$ Basic Partial Map \textsc{diagram} of 
\textbf{Structure Theorem} for $\hatPRa$ out of RCF\,1:

\bigskip
\begin{minipage} {\textwidth}
$
\xymatrix @+2em { 
D_{f}
\ar @/^1em/ [rrd]^{\widehat{f}}
\ar @<-0.5ex> [d]_{d_f}
\ar @{} [rd] |{\quad \pareq}  
  & 
\\
A 
\ar @{-^{>}} @<-0.5ex> [u]_{d^-_f}
\ar @{-^{>}} [rr]^{f}
  & & B
}
$  
\end{minipage}

\bigskip
In fact, with both ``structural'' $\hatPRa$ equations of the
\textsc{diagram,} we get from our $\hatPRa$ equation: 
\begin{align*}
& f \parcirc g \parcirc f 
      \pareq f \parcirc g \parcirc \widehat{f} \parcirc d_f^-
        \pareq  f \parcirc g \parcirc f \parcirc d_f \parcirc d_f^- \\
& \pareq f \parcirc d_f \parcirc d_f^- \quad \text{by equation above} \\
& \pareq \widehat{f} \parcirc d_f^- \pareq f: A \parto B \parto A \parto B
                                                       \quad\textbf{\qed}
\end{align*}  


\section{Choice within Classically Quantified Arithmetics}

\textbf{Define} Theory $\bfQ = \PRa+\forall\exists!$ as 
as Cartesian (!) PR extension of $\PRa$ by 
\emph{Quantification}---considered to give 
$\exists_b\,\ph(a,b): A \to 2$ and $\forall_b\,\ph(a,b): A \to 2$ 
as \emph{total} maps, this (intuitive) \emph{totality} 
\emph{formally} expressed by---\textbf{axiomatically}
maintained---\textbf{Cartesianness,} and by 
possibility---axi\-o\-ma\-ti\-cal\-ly forced either---of 
\textbf{map-definition} via (formal) \emph{unique existence} of 
\emph{values} to given \emph{arguments.} 

Formally, we \textbf{define} ``minimal classical'' (categorical) Theory 
$\bfQ$ by the following additional schemata over $\PRa:$ 

\medskip
$\bullet$ ``Quantified'' \emph{law of excluded middle:}  

\smallskip
\inference{ (\mr{no\text{-}mid}) }
{ $\ph = \ph(a,b): A \times B \to 2$ in $\bfQ$ }
{ $\bfQ \derives\ [\,\forall_b\,\ph(a,b)
                      \,\lor\,\exists_b\,\neg\,\ph(a,b)\,] 
                                        = \true_A(a): A \to 2$ }

\smallskip
$\bullet$ ``Discrete'' \emph{Map definition by unique existence:}
\inference{ (\forall\exists!) }
{ $\ph = \ph(a,b): A \times B \to 2$ 
         \emph{functional from} $A$ to $B,$ \ie \\
& $\bfQ \derives\ (\forall\,a \in A)\,(\exists!\,b \in B)\,\ph(a,b)$ \\
& \quad
    $[\,$\emph{Unique} existence is formalised as usual by a \\
& \quad
    Free-Variables implication between maps.$]$
}
{ $f = f_{\ph} = f_{\ph} (a): A \to B,$  
               \ in $\bfQ$ \emph{characterised} by \\ \\
& $\bfQ\ \derives\ [\,f_{\ph} (a) \doteq_B b\,] 
                         = \ph(a,b): A \times B \to 2$
}

Forgoing schema---including its \emph{uniqueness} clause---then gives,
for all $\bfQ$-maps $f,g: A \to B:$
\begin{itemize}
\item
\emph{argumentwise functionality:}
\begin{align*}
\bfQ\ \derives\ 
& \forall_a\,\exists!_b\ f(a) \doteq b,\ \text{in FV form, for}\ a \in A: \\ 
\bfQ\ \derives\ 
& \exists!_b\,[\,f(a) \doteq b\,]: A \to 2,\ a \in A\ \free
\end{align*}

\item
\emph{argumentwise definition} of \emph{map-equality:}
\begin{align*}
& \bfQ\ \derives\ f = g: A \to B\quad\text{\textbf{iff}} \\
& \bfQ\ \derives\ [\,f(a) \doteq_B g(a)\,]: A \to 2,\ a \in A\ \free 
    \quad (\emph{Equality Definability}) \\
& \textbf{iff} \quad \bfQ\ \derives\ \forall_a\,[\,f(a) \doteq_B g(a)\,]: 
                                                              \one \to 2
\end{align*} 

\end{itemize}

What we want to show for Theory $\bfQ$ is a (map-theoretical) 
\emph{local version} of the \emph{Axiom of Choice}, $\AC,$ 
necessarily here just---pointed---\emph{Countable Choice} $\ACC.$

  %




 

\medskip
\textbf{$\mu$-Inherit Lemma:} 
$\bfQ$ (and hence $\bfQ^{\wo}$) \textbf{inherit} $\hatPRa$'s $\mu$-operator: 
for $\PRa$-predicate $\ph = \ph(a,n): A \times \N \to 2,$
\begin{align*}
\bfQ\ \derives\ 
& \mu^{\PRa} \ph = \mu^{\PRa}\set{n\,|\,\ph(a,n)}: A \parto \N \\
& \pareq
  \mu^{\bfQ} \set{n\,|\,\ph(a,n)} \bydefeq  
  \begin{cases}
    \min\set{n'\,|\,\ph(a,n')}\ \myif\ \exists\,n\,\ph(a,n) \\
    \text{\emph{undefined}}\ \myif\ \forall\,n\,\neg\,\ph(a,n) 
  \end{cases} \\
& : A \parto \N
\end{align*}

\textbf{Proof:} 
Asserted partial-map \emph{equality}
$\bfQ\ \derives\ \mu^{\bfQ}\ph \pareq \mu^{\PRa}\ph: A \parto \N$
for $\PRa$-predicates $\ph = \ph(a,n): A \times \N \to 2$ is due to
the fact that the two $\mu$-recursive (partial) maps are 
\emph{compared}---in both directions---by suitable $\bfQ$-total 
maps with respect to their \emph{graphs,} as follows: 

\smallskip
Consider---within $\hatQ$---\textbf{defining} $\PRa$-diagram for \\
$\mu^{\PRa}\ph\,(a): A \parto \N,$ namely 
$$
\xymatrix @+2em{
\set{(a,n) \in A \times \N\,|\,\ph(a,n)}
\ar [d]^{a}_{\ell}
\ar [dr]^{\qquad\min\set{n' \leq n\,|\,\ph(a,n')}}
\\
A
\ar @{-^>} [r]^{\mu^{\PRa} \ph}
& \N 
}
$$
\emph{Partial} $\mu$-recursive map 
$\mu^{\bfQ}\ph = \mu^{\bfQ}\ph\,(a) \parto \N$
\textbf{defines}---``over'' $\bfQ$---an \emph{equal} partial map by
$$
\xymatrix @+2em{
\set{a \in A\,|\,\exists\,n\ \ph(a,n)}
\ar [d]^{\subseteq}
\ar [dr]^{\quad\min\set{n'\,|\,\ph(a,n')}}
\\
A
\ar @{-^>} [r]^{\mu^{\bfQ} \ph}
& \N 
}
$$
Partial-map-\emph{equality} \ $\mu^{\PRa}\ph \pareq \mu^{\bfQ}\ph$ 
\ in $\hatQ$ \,(``over $\bfQ$\,'') established by $\bfQ$-maps
\begin{align*}
& i :\,= \ell\ \circ\,\subseteq\,: 
    \set{A \times \N\,|\,\ph} 
         \to \set{a \in A\,|\,\exists\,n\ \ph(a,n)}\ \,\text{and} \\
& j :\,= (a,\mu^{\bfQ}\ph\,(a)): \set{a \in A\,|\,\exists\,n\ \ph(a,n)}
                                     \to \set{A \times \N\,|\,\ph}.
\end{align*}
Both $i,\,j$ \emph{total} $\bfQ$-maps, since
  $$\mu^{\bfQ}\ph\,(a) \bydefeq \min\set{n'\,|\,\ph(a,n')}: 
                  \set{a \in A\,|\,\exists\,n\ \ph(a,n)} \to \N$$
is $\bfQ$-\emph{total} \,\textbf{\qed}

\medskip\noindent
This \textbf{Lemma} gives

  \newpage

\medskip
\textbf{Middle-Inverse Countable-Choice Theorem} 
for Theory $\bfQ:$ 

\begin{itemize}

\item
Since $\bfQ$ \textbf{extends} $\PRa$ and $\hatPRa,$ it inherits
\emph{middle-inverse-property} from $\hatPRa.$ 
In particular for a \emph{pointed} $\PRa$-map
$f = f(a): A \to B$ (point $a_0: \one \to A$ given), $\bfQ$ \emph{inherits} 
earlier \emph{partial} $\mu^{\PRa}$ map
\begin{align*} 
f^- =\, 
& f^-(b) \bydefeq 
                 \cont_A(\mu^{\PRa} \set{n\,|\,f(\cont_A(n)) \doteq_B b}) \\
\pareq\, 
& \cont_A(\mu^{\bfQ} \set{n\,|\,f(\cont_A(n)) \doteq_B b}):
                                                  B \parto \N \to A
\end{align*}
as (partial) \emph{middle inverse.}

\smallskip
$[\,\cont_A (n): \N \to A$ retractive 
                \emph{count,} available via point $a_0\,]$ 

\item 
for a (pointed) $\PRa$-map $f: A \to B,$ the $\bfQ$-map 
$$
f' = f'(b): B \to A \\ 
\defeq
\begin{cases}
f^-(b)\ \,\myif\ f^-(b)\ \text{\emph{defined}, \ie} \\
\quad
  \myif\ \,\exists\,n\ f(\cont_A(n)) \doteq_B b \\
a_0\ \,\text{\emph{otherwise,}} \\
\quad
  \ie\ \myif\ \forall\,n\ f(\cont_A(n)) \neq_B b
\end{cases}
$$
is \textbf{definitionally} \emph{complemented} 
into a (``non-constructive'') \emph{total} $\bfQ$-map, $f': B \to A,$
\emph{middle-inverse} to given $f: A \to B$ in the sense
of schema $\ACCmi.$

\end{itemize}
  
\textbf{Comment:} We will \emph{not} rely on latter 
\emph{middle-inverse} $\bfQ$-choice map $f': B \to A$---which 
involves ``\,$\forall$\,'' and the Quantified law of excluded middle.









\medskip
$\forall$-\textbf{Elimination:} 
For our argument below we may drop \emph{formal} universal Quantor
``\,$\forall$\,'', Quantified \emph{law of excluded middle,} and 
\textbf{replace} schema $(\forall\exists!)$ above by 
\textbf{schema} of \emph{map definition by unique value-existence} 
\inference{ (\mr{FV}/\exists!) }
{ $\ph = \ph(a,b): A \times B \to 2$ \\ 
& \qquad  
    \emph{FV/$\exists!$ functional from} $A$ to $B,$ \ie \\
& $\bfQ \derives\ (\exists!\,b \in B)\,\ph(a,b): A \to 2,\ a \in A\ \free$
}
{ $f_{\ph} = f_{\ph} (a): A \to B,$ \\  
& \qquad
    this map $f_{\ph}$ in $\bfQ$ (``again'') characterised by \\
& $\bfQ\ \derives\ [\,f_{\ph} (a) \doteq_B b\,] 
                                = \ph(a,b): A \times B \to 2,$ \\
& \qquad\quad
    $a \in A,\ b \in B\ \free$
}  
\textbf{as well as} (canonical) 
     \emph{map definition via ``multivalued'' predicate,}
\inference{ (\mr{FV}/\exists) }
{ $\ph = \ph(a,b): A \times B \to 2$ in $\PRa,$ \\ 
& $\bfQ \derives\ (\exists\,b \in B)\,\ph(a,b): A \to 2,\ a \in A\ \free$
}
{ $\bfQ\ \derives\ f_{\ph} = f_{\ph} (a) 
         \defeq \mu^{\bfQ} \set{b \in B\,|\,\ph(a,b)}: A \to B$ 
                                                       \,(``total'') \\
& $ = \mu^{\PRa} \set{b \in B\,|\,\ph(a,b)}: A \parto B$ 
                               \,in\, $\hatPRa \bs{\prec} \bfQ$ 
                                                     \,(\textbf{subSystem})
} 
This $\bfQ$-\emph{map} $f_{\ph}$ is characterised within $\bfQ$ by 
\begin{align*}
\bfQ\ \derives\ 
& [\,f(a) \doteq_B b\,] = \ph(a,b): A \times B \to 2, \\
& a \in A,\ b \in B\ \free,\ \,\text{and \emph{value-minimality:}} \\ 
\bfQ\ \derives\ 
& [\,\ph(a,b) \impl f(a) \leq_B b\,]: A \times B \to 2
\end{align*}
in the order (canonically) inherited by $B$ (pointed) from that
of $\N,$ via retraction $\cont_B = \cont_B(n): \N \to B.$  
  
\medskip
So the critical properties of $\bfQ$ are those of its 
\emph{\textbf{existential}} Quantification. 
\emph{\textbf{first:}} this quantification yields
\emph{total predicates,} in the formal sense that it leads never out
of \emph{Cartesianness,} and \emph{\textbf{second:}} it allows---by sheer 
(established) \emph{formal existence} of ``values''---
\textbf{definition} of maps via (even ``infinite'') argument/value tables.

\section{Complexity Controlled Iteration Recalled}

Complexity Controlled Iteration---$\mathrm{CCI}_O$---is \emph{Iteration}
of a \emph{predecessor} (endo) step, decreasing \emph{Complexity}
of argument---Complexity measured in (a given) \emph{Ordinal} 
$O$---as long as complexity zero is not ``yet'' reached. \emph{Result} 
then is the argument \emph{reached,} with complexity zero. We choose 
here (\emph{axis} case) $O :\,= \N[\omega] \subset \omega^{\omega},$ 
the set of \emph{polynomial coefficient} strings (no trailing zeros). 

It is highly plausible, and a Theorem in $\PA$---at least
in $\PA+\text{\emph{well-order} of}\ \omega^{\omega}$ --, that such 
$\mathrm{CCI}$'s \emph{terminate,} on each initial argument 
given. So our first step in direction of \emph{Terminating} 
Recursiveness---\emph{strengthening} $\PRa$---can (and will) 
be formalisation first of the concept $\mrCCI$ of 
\emph{Complexity Controlled Iteration} (``over'' $O :\,= \N[\omega]$) 
and---second---introduction of \emph{axiom} schema for 
conceiving \emph{weakest} Theory $\piR$ (strengthening $\PRa$ and)
admitting \emph{termination} of all these $\mrCCI$'s.

We attempt to formalise wanted Theory within the \emph{partial-map} 
framework of theory $\hatPRa \bs{\sqsupset} \PRa,$ 
which is a \emph{definitional, conservative} extension of
Theory $\PRa.$ It contains (Cartesian) $\PRa$ \emph{embedded} as  
a \emph{monoidal PR subCategory.}

\smallskip
\textbf{Definition:} For ``\emph{Ordinal}'' $\N[\omega],$ 
schema $(\mrCCI)$ below (quote from part RCF\,2) is to \textbf{define} 
a \emph{Complexity Controlled Iteration}---$\mrCCI$---with 
\emph{complexity values} in $\N[\omega],$ as a (formally) 
\emph{partial} map $\wh[c>0\,|\,p]: A \parto A,$ definition based 
on suitable data $c$ (complexity) and $p$ (predecessor step) 
as follows, within any theory $\bfS$ strengthening $\PRa:$ 
\inference{ (\mrCCI) }
{ $c: A \to O$ in $\bfS,$\quad \emph{complexity}, \\
& $p: A \to A$ $\bfS$-endo,\quad \emph{predecessor step,} \\
& $\bfS \derives\ c(a) > 0_O \implies p\ c(a) < c(a): A \to 2$
                                       \quad (\emph{Desc}) \\ 
& \qquad 
    \emph{strict descent above complexity zero}, \\
& $\bfS \derives\ c(a) \doteq 0_O \implies p(a) \doteq_A a: A \to 2$
                                       \quad (\emph{Stat}) \\ 
& \qquad 
    \emph{stationarity at complexity zero}
}
{ $\wh[\,c>0\,|\,p\,]: A \parto A$ in $\hatS$ (partial) map: 
}
$\wh[\,c>0\,|\,p\,]: A \parto A$ realises the $\mrCCI$ (as a $\while$ loop).
As a partial map it is given \emph{defined arguments enumeration}
\begin{align*}
& d_{\wh[\,c>0\,|\,p\,]}\,(a,n) \bydefeq a: \\
& D_{\wh[\,c>0\,|\,p\,]} 
    = \set{(a,n) \in A \times \N\,|\,c\ p^n(a) \doteq 0_O} 
            \xto{\subseteq} A \times \N \xto{\ell} A,
\end{align*}
and (calculation) \emph{rule}
  $$\widehat{\wh} [\,c>0\,|\,p\,] = \widehat{\wh} [\,c>0\,|\,p\,]\,(a,n) 
                               \defeq p^n(a): D_{\wh[\,c>0\,|\,p\,]} \to A.$$
\textbf{Comment:} Essential ``ingredient'' for above iteration 
$\wh[\,c>0\,|\,p\,]$ is its (formally) \emph{partial} \emph{termination-index} 
\begin{align*}
& \mu[\,c>0\,|\,p\,] = \mu[\,c>0\,|\,p\,] (a): A \parto \N
                                              \ \,\text{given as} \\
& \mu[\,c>0\,|\,p\,] (a) \defeq \mu\set{n\,|\,c\ p^n(a) \doteq 0}: 
                                                       A \parto \N,
\end{align*}
and as such characterised---as \emph{partial} map: within $\hatS$---by 
\begin{align*}
& \hatS \derives\ c\ p^\S(a,\mu[\,c>0\,|\,p\,]\,(a)) \doteq 0,
                                                      \ \text{and} \\ 
& \mu[\,c>0\,|\,p\,]\,(a): A \parto \N
           \ \text{(argumentwise) \emph{minimal} in this regard.}
\end{align*}
\text{Partial map}
  $$ (\id_A,\mu[\,c>0\,|\,p\,]) 
        = (A \times \mu[\,c>0\,|\,p\,]) \parcirc \Delta:
                                       A \to A^2 \parto A \times \N$$
is just the---pointwise minimised (``canonical'')---\emph{opposite} 
partial map
  $$d^- = d_{\wh[\,c>0\,|\,p\,]}^-: A \parto D_{\wh[\,c>0\,|\,p\,]},$$ 
opposite to $d = d_{\wh[\,c\,|\,p\,]}.$ As \emph{opposite,} this 
$d^-$ has \emph{partial section} property 
$d \parcirc d^-\ \parinc\ \id_A$ within $\hatS,$ \emph{maximally.}






\section{Cartesian Code Self-Evaluation ``inside'' $\bfQ^{\wo}$} 

We \emph{question} here---on \emph{consistency}---Theory 
$\bfQ^{\wo} \defeq \bfQ+\wo(\omega^{\omega})$ \\ 
$\bydefeq \PRa+\forall\exists!+\wo(\omega^{\omega})$ of 
\emph{classically Quantified} Arithmetic with well-ordered 
$\omega^{\omega},$ subsystem of \textbf{set} theory $\ZF.$

\medskip
We attempt to exhibit a \emph{Cartesian} \textbf{subsystem} 
$\dotpiR$ of $\bfQ^{\wo}$ which admits ``total'' (\ie Cartesian)
\emph{self-evaluation} $\doteps = \doteps(u,a): \mrdotpiR \times \X \to \X.$
As a \textbf{consequence,} $\bfQ^{\wo}$ will turn out to be 
\emph{inconsistent.}

\smallskip
We first form the \emph{monoidal closure} 
$\piR+(d^-_{\wh[\,c>0\,|\,p\,]}),$ within $\hatpiR,$ of $\piR$ under
all (formally partial) $\mu$-recursive maps \emph{opposite,}
(canonically) \emph{middle-inverse,} to the (PR) 
\emph{defined-arguments} enumerations
  $$d_{\wh[\,c>0\,|\,p\,]}: D_{\wh[\,c>0\,|\,p\,]} 
     = \set{(a,n)\,|\,c\,p^n \doteq 0 \in \N[\omega]} \ovs{a} A$$
of all $\mrCCI$'s given by (PR) \emph{complexity} $c$ and 
$\N[\omega]$-descending (PR) step $p: A \to  A.$

Interpreted in \textbf{frame} $\bfQ^{\wo},$ these $d^-_{\wh[\,c>0\,|\,p\,]}$
become \emph{total,} since
\begin{align*}
\bfQ^{\wo}\ \derives\ 
& \exists\,n\,[\,c\,p^n(a) \doteq 0\,]: A \to 2,\ (a \in A\ \free), 
                                                     \ \text{hence} \\
\bfQ^{\wo}\ \derives\ 
& d^-_{\wh[\,c>0\,|\,p\,]} (a) \bydefeq 
              (a,\mu\set{n\,|\,c\,p^n(a) \doteq 0\,}) \\
& \bydefeq (a,\mu^{\PRa} \set{n\,|\,c\,p^n(a) \doteq 0\,}) \\
& = (a,\mu^{\bfQ^{\wo}} \set{n\,|\,c\,p^n(a) \doteq 0\,}):
                                     \ (\text{$\mu$-\textbf{inherit}}) \\
& \quad
    A \parto D_{\wh[\,c>0\,|\,p\,]}\ \,\text{\emph{total}}: \\
& \quad
    A \to D_{\wh[\,c>0\,|\,p\,]} \subset A \times \N
\end{align*}
and they are, again within $\bfQ^{\wo},$ \emph{sections} to their 
  $$d_{\wh[\,c>0\,|\,p\,]}: D_{\wh[\,c>0\,|\,p\,]} \to A,$$
(not only \emph{partial} sections, \cf the above). 

\smallskip
\emph{Categorically,} this totality means that the \NAME{Godement}
equations hold, ``even'' when these \emph{additional} $\hatPRa$-maps 
$d^-$ are involved.

\smallskip
So the following Theory $\dotpiR,$ strengthening of monoidal theory \\
$\piR+(d^-_{\wh[\,c>0\,|\,p\,]}) \bs{\sqsubset} \hatpiR$ by the following 
two \textbf{axioms}
(``schemata'') is \\
---as a \textbf{subsystem} of $\bfQ^{\wo}$---consistent 
relative to $\bfQ^{\wo}:$
\inference{ (Gode) }
{ $f: C \parto A,\ g: C \parto B$ in $\piR+(d^-_{\wh[\,c>0\,|\,p\,]})$ }
{ For ``induced'' $(f,g) \bydefeq (f \times g) \parcirc \Delta_C:$ \\
& \qquad\qquad 
    $C \to A \times A \parto B \times B:$ \\
& $\dotpiR\ \derives\ \ell \circ (f,g) = f: 
                                  C \to A \times B \ovs{\ell} A$ and \\
& $\dotpiR\ \derives\ r \circ (f,g) = g: C \to A \times B \ovs{r} B$
}
$[\,$We could add as an axiom \emph{section} property
$d_{\wh}\,\circ\,d^-_{\wh[\,c>0\,|\,p\,]} = \id_A,$ given wihin $\bfQ^{\wo}$
as well, but forcing \emph{Cartesianness}---standing for totality of
all of $\dotpiR$ maps---will be sufficient for our argument$\,]$
 

\medskip
Since evaluation $\eps = \eps(u,a): \mrPRa \times \X \parto \X$
is \textbf{defined} as a $\mrCCI$ within $\hatPRa,$ 
\begin{align*}
& d^-_{\eps} = d^-_{\eps} (u,a)
    \bydefeq ((u,a),\mu\set{n\,|\,\cPR\,e^n(u,a) \doteq 0}): \\
& \qquad
    \mrPRa \times \X \to D_{\eps} \subset (\mrPRa \times \X) \times \N
\end{align*}
is in fact a $\dotpiR$-map (considered \emph{total} as such).

\smallskip
Is it possible to extend this $\PRa$-evaluation 
$\eps(u,a): \mrPRa \times \X \parto \X$---defined as a $\hatPRa$-map---into
a \emph{code self-evaluation} 
  $$\doteps = \doteps(u,a): \mrdotpiR \times \X \to \X\,?$$
For this end, let us treat the \emph{additional} maps 
  $$d^- = d^-_{\wh[\,c>0\,|\,p\,]}: 
                    A \to D = D_{\wh[\,c>0\,|\,p\,]} \subset A \times \N$$
as ``basic'' with respect to the evaluating $\mrCCI$ to be constructed: 
\begin{align*}
& \dot{c}\,(u,a) = c_{\doteps}\,(u,a): 
                  \dotpiR \times \X \ovs{u} \dotpiR \to \N[\omega] \\
& \qquad
    \text{is PR \textbf{defined} from}\ \cPR: 
                   \mrPRa \times \X \to \mrPRa \to \N[\omega] \\
& \qquad
    \text{by adding the clause} \\
& \wh = \wh[\,c>0\,|\,p\,]\ \,\mrCCI 
          \implies \dot{c}\,(\code{d^-_{\wh}}) \doteq 1.
\end{align*}
Extended evaluation \emph{step} 
\begin{align*}
& \dot{e} = \dot{e}\,(u,a): \dotpiR \times \X \to \dotpiR \times \X \\
& \qquad
    \text{then is PR \textbf{defined} from}\ 
                       e: \mrPRa \times \mrPRa \to \mrPRa \times \mrPRa \\
& \text{by addition of (\emph{Objectivity}) clause} \\
& \wh = \wh[\,c>0\,|\,p\,]\ \,\mrCCI 
    \implies 
      \dot{e}\,(\code{d^-_{\wh}},a) \defeq (\code{\id}, d^-_{\wh}\,(a)) \\ 
& [\ \in \mrdotpiR \times (A \times \N) \subset \mrdotpiR \times \X\ ]
\end{align*}
(Self-) evaluation $\doteps = \doteps\,(u,a): \mrdotpiR \times \X \to \X$ 
then is \textbf{defined}---within Cartesian theory $\dotpiR$ itself----by 
$\mrCCI$ (!)
\begin{align*}
& \doteps = \doteps\,(u,a) \defeq r \circ \wh[\,\dot{c}\,|\,\dot{e}\,]\,(u,a):
    \mrdotpiR \times \X \to \mrdotpiR \times \X \ovs{r} \X \\
& \qquad
    \text{with $\hatPRa$ middle-inverse} \\
& d^-_{\doteps}\,(u,a) = d^-_{\wh[\,\dot{c}>0\,|\,\dot{e}\,]}\,(u,a): 
      \dotpiR \times \X \to D_{\doteps} 
                 = \set{((u,a),n)\,|\,\dot{c}\,\dot{e}^n(u,a) \doteq 0}
\end{align*}
\text{a $\hatPRa$ (partial) map, more: a $\dotpiR$ \emph{map.}}  
By \textbf{Structure Theorem} for $\PRa,$ we have for any (partial)
$\hatPRa$-map $f = \an{(d_f,\widehat{f}): D_f \to A \times B}: A \parto B:$
  $$\hatPRa\ \derives\ 
      f \pareq \widehat{f} \parcirc d^-_f: A \parto D_f \to B.$$
So for $\mrCCI$'s: 
\begin{align*}
\pi\widehat{\dot{\mathbf{R}}}\ \derives\ 
& \wh[\,c>0\,|\,p\,] = \widehat{\wh}\,\parcirc\,d^-_{\wh}: 
                                          A \parto D_{\wh} \to A \\
& \text{is \emph{represented} as a $\dotpiR$-map}\ \dot{\mr{w}}\mr{h}: \\
\pi\widehat{\dot{\mathbf{R}}}\ \derives\ 
& \wh[\,c>0\,|\,p\,] \pareq \dot{\mr{w}}\mr{h} 
           \defeq \widehat{\wh} \circ d^-_{\wh}:  A \to D_{\wh} \to A.
\end{align*}
As a special case then 
  $$\dotpiR\ \derives\ 
      \doteps = \doteps\,(u,a) = \widehat{\doteps} \circ d^-_{\doteps}\,(u,a):
                  \dotpiR \times \X \to D_{\eps} \ovs{\widehat{\doteps}} \X$$
is a (is represented as) map in $\dotpiR$ itself. It
constitutes a (code) \emph{self-evaluation} for Theory $\dotpiR$ 
since---only further property needed---it is \emph{Objective} as 
an evaluation, will say
\inference{ (\mr{Ob}_{\doteps}) }
{ $f = f(a): A \to B$ a $\dotpiR$-map }
{ $\dotpiR\ \derives\ \doteps(\code{f},a) = f(a): A \to B$ }

\textbf{Proof:} Objectivity 
$\PRa\ \derives\ \eps\,(\code{f},a) = f(a): A \to B$
of ``fundamental'' evaluation $\eps: \mrPRa \times \X \parto \X$ 
has been shown in RCF\,2, by external \ultxt{PR} on $\ultxt{depth}\,[f]$ 
which relies on \textbf{Peano Induction} (on the itertion counter $n$)
in case $(\code{f},a) \in \mrPRa \times \X$ of form of an iterated:
$(\code{f},a) = (\code{g^{\S}},\an{b;n}).$

\smallskip
$[\,$Free-Variables Peano Induction is available 
                            in $\PRa$ and strengthenings$\,]$

\smallskip
Same \ultxt{PR} \ultxt{argument} works in present case of self-evaluation 
  $$\doteps\,(u,a) = r \circ \wh[\dot{c}>0\,|\,\dot{e}\,]: 
      \mrdotpiR \times \X \to \mrdotpiR \times \X \xto{r} \X.$$
The reason is that the evaluation clause for the additional
maps is given as an \emph{Objectivity} instance:
\begin{align*}
& \dot{e}\,(\code{d^-_{\wh[\,c>0\,|\,p\,]}},a) 
                           \bydefeq (\code{\id},d^-_{\wh}\,(a)), \\
& \dotpiR\ \derives\ 
    \doteps\,(\code{d^-_{\wh[\,c>0\,|\,p\,]}},a) = d^-_{\wh}\,(a): \\
& \qquad\quad
    \X \supset A \to D_{\wh} \subset A \times \N \subset \X
\end{align*}  
So Objectivity is preserved by extension of evaluation
$\eps$ to $\doteps: \mrdotpiR \times \X \to \X$ \,\textbf{\qed}
 
\medskip
But (Objective) \emph{code self-evaluation} of any \emph{Cartesian} PR 
Theory $\T$ renders $\T$ \textbf{inconsistent,} as we will show in 
detail---final section---by the ``appropriate'' \emph{diagonal} argument.

Since self-evaluating theory $\bfQ^{\wo}$ is an extension of 
inconsistent $\dotpiR,$ $\bfQ^{\wo}$ itself turns out to be
inconsistent. So in particular Peano Arithmetic $\PA+\wo(\omega^{\omega})$
with the lexicographical order on $\omega^{\omega} = \N^* \supset \N[\omega]$ 
a \emph{well-order}, as well as \textbf{set} theory $\ZF$ are
shown to be inconsistent.

\section{Liar via Code Self-Evaluation}

Any \emph{Code Self-Evaluation} family
  $$\eps = \eps_{A,B}\,(u,a): \cds{A,B}_\T \times A \to B$$
of a \emph{Cartesian} PR Theory $\T$ within Theory $\T$ \emph{itself,}
which is \emph{Objective} as (self-) evaluation---see above---establishes 
a \textbf{contradiction} within $\T,$ by the (``usual'') diagonal 
argument below: formalisation of ``Antinomie Richard''.


\smallskip
Remains to develop that \emph{diagonal} argument ``against''
(consistent) \emph{code self-evaluation} for Theory $\T$
in general---skip, if you are used to such \emph{diagonal argument}
--, same argument as in RCF\,3: \emph{Map-Code Interpretation via Closure.} 

In presence of such (Objective) \emph{self-evaluation}
family $\eps$ \textbf{define} (anti) \emph{diagonal} 
$d: \N \to 2$ within general \emph{Cartesian} Arithmetical theory $\T:$ 

\smallskip
\qquad
  $d \defeq \neg \circ \eps_{\N,2}\,\,\circ\,(\#,\id_\N): 
      \N \lto \cds{\N,2}_{\T} \times \N \lto 2 \overset{\neg} {\lto} 2,$

\smallskip\noindent
with $\# = \#(n): \N \overset{\iso} {\lto} \cds{\N,2}_{\T}$ the 
---isomorphic---PR \emph{count} of all (internal) predicate codes,
of Theory $\T.$ 

\smallskip
As expected in such \emph{diagonal argument,} we \emph{substitute} --
within Theory $\T$---the counting index 
  $$q \defeq \#^{-1} (\code{d}) = \#^{-1} \circ \code{d}: 
                \one \to \cds{\N,2}_{\T} \overset{\iso} {\lto} \N$$
of $d$'s \emph{code,} into $\T$-map $d: \N \to 2$ itself, 
and get a ``liar'' map $\liar: \one \to 2,$ namely
\begin{align*}
\T\ \derives\ 
& \liar \defeq d\,\circ\,q: \one \to \N \to 2 \\
& \bydefeq d\,\circ\,\#^{-1}\,\circ \code{d} \\
& \bydefeq \neg\,\circ\,\eps_{\N,2}\,\circ\,(\#\,,\,\id_\N) 
               \,\circ\,\#^{-1}\,\circ \code{d} \\
& = \neg\,\circ\,\eps_{\N,2}
        \,\circ\,(\code{d}\,,\,\#^{-1}\,\circ \code{d}) 
                                 \qquad (\T\ \text{\emph{Cartesian}}) \\
& \bydefeq \neg\,\circ\,\eps_{\N,2}\,(\code{d}\,,\,q) \\
& = \neg\,\circ\,d(q)
             \qquad (\emph{Objectivity of}\ \eps) \\
& \neg\,\circ\,d\,\circ\,q \bydefeq \neg\ \liar: \one \to 2 \to 2,
\end{align*}
a \emph{contradiction,} whence

\smallskip
\textbf{Conclusion} (again): Code-self-evaluating Theory $\dotpiR$ 
is \textbf{inconsistent} and so are all of its \textbf{extensions,}
in particular ``minimal'' Quantified Arithhmetical Theory
$\bfQ^{\wo} = \bfQ+\wo(\omega^{\omega})$ with $\forall\,\exists!$ 
\textbf{definition} of maps out of (binary, PR) predicates as well 
as \textbf{its} \textbf{extensions} such as $\PA+\wo(\omega^{\omega}),$ 
and extension $\ZF$ of the latter theory.

\smallskip
$[\,$Without well-order of ``one of the first'' countable 
``\emph{Ordinals}'', namely of $\omega^{\omega},$ (countable) 
well-order which is expressible within the language of first-order 
\textbf{set}-theory $\oneZF$ and already within that of $\PA,$ the 
theory of (countable) \emph{Ordinals} would be rather poor$\,]$

\section{Discussion} 

$\bullet$
Our \textbf{inconsistency} argument applies to 
\textbf{Peano-Arithmetic,} if this theory is presented as
predicate calculus (``full quantification'') for description of
\emph{Algebra} \& \emph{Order} on $\N,$ plus induction schema 
$\mr{P5,}$ \ie if $\PA$ is conceived as $\PR\forall\exists!:$
PR \emph{infinity} plus ``full'' (classical) predicate calculus,
with ``set theoretical'' possibility of \emph{map-definition,} 
see above.

But following \NAME{Lawvere}---and \NAME{Goodstein}---Algebra 
\& Order can be expressed by equations, in particular 
by use of \emph{truncated subtraction}
$[\,m \dotminus n\,]: \N \times \N \to \N,$ which yields 
order and equality predicates on $\N \times \N,$ as well as 
(\emph{constructive}) ``existence'' of $b$ such that 
$a+b = c,$ namely $b :\,= c \dotminus a$ for ``given'' 
$a,c \in \N$ satisfying $a \leq c,$ see the \NAME{Wikipedia}-article 
on \emph{Peano Arithmetic,} and the fact that Peano induction 
$\mr{P5}$ can be expressed equationally, within (categorical) 
Free-Variables Calculus. 

\smallskip
$[\,$In $\PRa,$ induction axiom $\mr{P5}$ is a consequence of 
\emph{uniqueness} of maps \emph{defined} by the \emph{full schema} 
of Primitive Recursion.$]$ 
 
\smallskip
But in Free-Variables setting---my guess---the PR schema, 
in form of (special one-fold successor case of) 
\emph{iteration schema} of \NAME{Eilenberg} \& \NAME{Elgot,} 
plus \NAME{Freyd}'s \emph{uniqueness schema} for the 
\emph{initialised iterated,} are \textbf{needed} for (unique) 
\textbf{definition} of the \emph{more complex} PR maps such as 
\emph{exponentiation,} \emph{faculty} \etc which are classically 
obtained from addition, multiplication, and order by use of 
(formal) \emph{existence.}   

\medskip
$\bullet$
On the ``constructive'' side, Free-Variables categorical Primitive 
Recursion Theory $\PRa$ above, strengthens into theorie(s) 
$\piOR$ ($O$ an Ordinal extending $\N[\omega]$) of
\emph{on-terminating} (\emph{not:} ``retractively'' terminating) 
\emph{Complexity-Controlled Iterations} with complexity
measured in $O:$ These theories ``just'' exclude 
\emph{infinite descending chains} in ``their'' Ordinal $O,$
and seem therefore to be almost as consistent as \emph{basic} 
Theory $\PRa$ (conservative extension of \emph{fundamental} 
Theory $\PR.$) Theories $\piOR$ are---``on the other
hand''---strong enough to \ultxt{derive} their own 
(Free-Variable) \emph{Consistency formulae,} see part RCF\,2 
mentioned above.

\medskip
$\bullet$
\textbf{Question:} Does our inconsistency argument equally apply 
to \emph{Arithmetical} \emph{first order Elementary Theory of Topoi} 
$\oneETT\N$ (Topoi with NNO) in place of Theory $\bfQ\,?$ 
As far as I can see, our argument could possibly be adapted to this 
case. Theory $\oneETT\N$ has \emph{two} \emph{truth}-Objects, one
arithmetical, $2 \bydefeq \set{\N\,|\,\,< 2} = \one \oplus \one,$ 
inherited from its \textbf{subSystem} $\PRa,$ as well as its
\emph{genuine,} intuitionistic \emph{subobject classifier} 
\xymatrix{
\one
\ar @<+0.5ex> [r]^{\false}
\ar @<-0.5ex> [r]_{\true}
& \Omega
} for its
``specific'' logic, in particular ``receiving'' (intuitionistic)
Quantifier $\exists.$  
$$
\xymatrix @+2em{
& \one
  \ar [d]_{0}
  \ar [drr]^{\false}
\\
2\ \text{is embedded into}\ \Omega\ \text{via} 
& 2 = \one \oplus \one
  \ar [rr]^{(\false\,|\,\true)}
  & & \Omega
\\
  & \one
    \ar [u]_{s\,0}
    \ar [urr]_{\true}
}
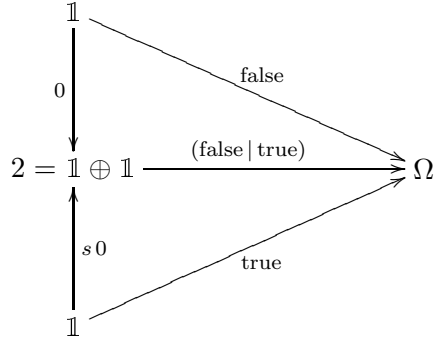
$$
$\oneETT\N$ admits schema
\inference{}
{ $\ph = \ph(a,n): A \times \N \to 2$ \\
& \emph{arithmetical} predicate, in $\PRa$ 
}
{ $\exists_r\ph = \exists\,n\,\ph(a,n): A \to \Omega$ 
    a (``total'') $\oneETT\N$-\emph{map}  \\
& ---``\,$\exists$ fits (already) in \emph{Cartesian} \textbf{frame}''--- \\
& + universal properties characterising map $\exists_r\ph: A \to \Omega$ \\ 
& within $\oneETT.$
}
My \textbf{guess} is that $\oneETT\N$ further admits \textbf{schema} of
$\exists$-\emph{dominated} \\
$\mu$-\emph{recursion}
\inference{ (\dot{\mu}) }
{ $\ph(a,n): A \times \N: A \times \N \to 2$ in $\PRa$ \\
& $\oneETT\N\ \derives\ \exists_r\ph = \true_A: A \to \Omega$
}
{ Formally \emph{partial} $\PRa$-map 
           $\mu_r\ph = \mu\set{n\,|\,\ph(a,n)}: A \parto \N$ \\
& \quad
    ``total'', \ie represented by a $\oneETT\N$ \emph{map} 
                                            $\dot{\mu}_r\ph: A \to \N$ \\
& $[\,$ \textbf{and} 
     $\oneETT\N\ \derives\ \ph(a,\dot{\mu}_r\ph\,(a)) = \true_A(a): A \to 2,$ \\
& $\dot{\mu}_r\ph: A \to 2$ (argumentwise) minimal in this regard$\,]$
}
``Latter \textbf{instance} of (overall) \emph{defined} 
$\dot{\mu}_n\ph(a,n): A \to 2,$ fits into (given) Cartesian frame 
of $\oneETT\N.$''

\smallskip
We saw above that we do not need formal \emph{universal} Quantor 
``\,$\forall$\,'', and in particular \emph{not} \emph{Booleanness}
of Quantification---we could drop schema of \emph{Excluded Middle.}

\smallskip
So, \textbf{if} $\oneETT\N$ should admit latter schema $(\dot{\mu}),$ 
of $\exists$-\emph{dominated totality} of $\mu_r\ph: A \to \N,$ 
\textbf{then} our \textbf{inconsistency} argument \emph{would} apply
to \emph{first order arithmetical Theory} $\oneETT+\wo(\omega^{\omega})$ 
of \emph{Elementary Topoi,} with \emph{lexicographical Order on 
$\omega^{\omega}$ a well-order.}

\smallskip
\textbf{If} so, then the final \textbf{question} is: Do \emph{real-life}
\textbf{Topoi,} \ie interesting Topoi of \textbf{sheaves,} have an NNO\,? 
  









\bigskip

\bigskip


  \noindent Address of the author: \\
  \NAME{M. Pfender}                       \hfill D-10623 Berlin \\
  Institut f\"ur Mathematik                              \\
  Technische Universit\"at Berlin         \hfill pfender@math.TU-Berlin.DE\\
  

\vfill

             \end{document}